\documentclass[a4paper,11pt]{amsart}
\usepackage{amsmath,amsthm,amssymb}
\usepackage[mathscr]{eucal}
 \usepackage{cite}
\usepackage{upgreek}
\usepackage{enumerate}
\usepackage[bookmarks=false]{hyperref}
\usepackage{enumitem}
\usepackage{amsmath}
\usepackage{mathrsfs}
\usepackage{blindtext}
\usepackage{scrextend}
\usepackage{enumitem}
\addtokomafont{labelinglabel}{\sffamily}
\usepackage{color}
\usepackage[at]{easylist}

\makeatletter
\newcommand\footnoteref[1]{\protected@xdef\@thefnmark{\ref{#1}}\@footnotemark}
\makeatother

\numberwithin{equation}{section}

\theoremstyle{plain}
\newtheorem*{theorem-non}{Main theorem}

\newtheorem{theorem}{Theorem}[section]
\newtheorem{claim}{Claim}
\newtheorem{case}{Case}
\newtheorem{lemma}[theorem]{Lemma}

\theoremstyle{definition}
\newtheorem{definition}[theorem]{Definition}
\newtheorem*{definition*}{Definition}

\newtheorem{notation}[theorem]{Notation}
\newtheorem{remark}[theorem]{Remark}

\begin{document}

\title[On sofic approximations of $\mathbb F_2\times\mathbb F_2$]
{On sofic approximations of $\mathbb F_2\times\mathbb F_2$}

\author[A. Ioana]{Adrian Ioana}
\address{Department of Mathematics, University of California San Diego, 9500 Gilman Drive, La Jolla, CA 92093, USA}
\email{aioana@ucsd.edu}

\thanks{The author was supported in part by NSF Career Grant DMS \#1253402 and NSF FRG Grant \#1854074.}
\begin{abstract}   We construct a sofic approximation of $\mathbb F_2\times\mathbb F_2$ that is not essentially a ``branched cover" of a sofic approximation  by homomorphisms. 
This answers a question of L. Bowen.

\end{abstract}

\maketitle

\section{Introduction} 
A countable group $\Gamma$ is called sofic if it admits  a sequence of almost actions on finite sets which are asymptotically free. To make this precise, endow the symmetric group $\text{Sym}(X)$ of any finite set $X$ with the normalized Hamming distance:
$\text{d}_{\text{H}}(\sigma,\tau)=|X|^{-1}\cdot |\{x\in X\mid\sigma(x)\not=\tau(x)\}|$. 
\begin{definition}
A sequence of maps  $\sigma_n:\Gamma\rightarrow\text{Sym}(X_n)$, for finite sets $X_n$, is called an {\it asymptotic homomorphism} if  $\lim\limits_{n\rightarrow\infty}\text{d}_{\text{H}}(\sigma_n(g)\sigma_n(h),\sigma_n(gh))=0$, for all $g,h\in\Gamma$. An asymptotic homomorphism $\sigma_n:\Gamma\rightarrow\text{Sym}(X_n)$ is called a {\it sofic approximation} of $\Gamma$ if it satisfies that $\lim\limits_{n\rightarrow\infty}\text{d}_{\text{H}}(\sigma_n(g),\text{Id}_X)=1$, for all $g\in\Gamma\setminus\{e\}$. The group $\Gamma$ is called {\it sofic} if it  has a sofic approximation. 
\end{definition}

In recent years, the study of sofic groups has received a lot of attention.
It is now understood that soficity has a number of important consequences (see, e.g., \cite{Bo18,Th18}). This is particularly interesting because sofic groups form a broad class, which includes all amenable and all residually finite groups.
Moreover, it is a longstanding open problem whether every countable group is sofic.

This note is motivated by the problem of classifying the sofic approximations of a given sofic group. 
For amenable groups $\Gamma$, a satisfactory classification of sofic approximations was found in \cite{EL11}:  any sofic approximation of $\Gamma$ is equivalent to one constructed from a disjoint union of F{\o}lner sets. Here, we say that two sofic approximations $\sigma_n:\Gamma\rightarrow\text{Sym}(X_n)$ and $\tau_n:\Gamma\rightarrow\text{Sym}(X_n)$ are {\it equivalent} if $\lim\limits_{n\rightarrow\infty}\text{d}_{\text{H}}(\sigma_n(g),\tau_n(g))=0$, for all $g\in\Gamma$ \cite{Bo17}. 
If $\Gamma$ is a residually finite group, then it admits a sofic approximation $\sigma_n:\Gamma\rightarrow\text{Sym}(X_n)$, where each $\sigma_n$ is a homomorphism.   Conversely, given a residually finite group $\Gamma$, one would ideally like to show that any sofic approximation of $\Gamma$ is equivalent to one consisting of homomorphisms, and thus arises from the finite quotients of $\Gamma$.
In this case,  $\Gamma$ is called {\it weakly stable} \cite{AP15} \footnote{This is a weakening of the notion of {\it stability in permutations} (or P-{\it stability}), requiring that any asymptotic homomorphism is equivalent to one given by homomorphisms.  For a survey of recent progress on stability, see  \cite{Io19b}.}. 
The class of weakly stable groups includes all residually finite  amenable groups \cite{AP15} and the free groups. 
As shown in  \cite{LLM19}, surface groups satisfy a flexible variant of weak stability.
On the other hand, we proved in \cite[Theorem D]{Io19b} that the product of two non-abelian free groups is not weakly stable. Consequently,  $\mathbb F_2\times\mathbb F_2$ admits a sofic approximation which does not essentially come from a sequence of homomorphisms. 

Our goal here is to strengthen this result and show the failure of a more general possible classification of sofic approximations of $\mathbb F_2\times\mathbb F_2$ proposed by L. Bowen. This is formulated using the following:

\begin{definition}\label{branch} 
Let $\sigma_n:\Gamma\rightarrow\text{Sym}(X_n)$ and $\tau_n:\Gamma\rightarrow\text{Sym}(Y_n)$ be asymptotic homomorphisms of a countable group $\Gamma$. We say that $(\sigma_n)$ is a {\it branched covering} of $(\tau_n)$ if there are onto maps $\theta_n:X_n\rightarrow Y_n$ such that $\theta_n\circ\sigma_n(g)=\tau_n(g)\circ\theta_n$, for all $g\in\Gamma$, and $\theta_n$ is a $d_n$-to-one, for some $d_n\in\mathbb N$.

\end{definition}
\begin{remark} \label{remark}
 Assume the setting from Definition \ref{branch}.  Then $\text{d}_{\text{H}}(\sigma_n(g),\text{Id}_{X_n})\geq\text{d}_{\text{H}}(\tau_n(g),\text{Id}_{Y_n})$, for every $g\in\Gamma$. 
Thus, if $(\tau_n)$ is a sofic approximation of $\Gamma$, then so is $(\sigma_n)$. 
 The branched covering construction therefore provides a way of producing new sofic approximations from  old ones.
 
We also  remark that any branched covering $(\sigma_n)$ of $(\tau_n)$  arises 
from a sequence of 
``almost cocycles" for $(\tau_n)$.
Indeed, let $Z_n=\{1,2,...,d_n\}$ and identify $X_n=Y_n\times Z_n$ so that $\theta_n:X_n\rightarrow Y_n$ is the projection map. Then $\sigma_n(g)(y,z)=(\tau_n(g)y,c_n(g,y)z)$, where $c_n:\Gamma\times Y_n\rightarrow\text{Sym}(Z_n)$ is a map satisfying $\lim\limits_{n\rightarrow\infty}|Y_n|^{-1}\cdot |\{y\in Y_n\mid c_n(gh,y)\not=c_n(g,\tau_n(h)y)c_n(h,y)\}|=0$, for all $g, h\in\Gamma$.

\end{remark}

At an Oberwolfach workshop in May 2011, Bowen asked  (see \cite[page 1463, Question 7]{OWR11}) if any sofic approximation $\sigma_n:\Gamma\rightarrow\text{Sym}(X_n)$ of $\Gamma=\mathbb F_2\times\mathbb F_2$ is essentially a branched covering of some sofic approximation $\tau_n:\Gamma\rightarrow\text{Sym}(Y_n)$ by homomorphisms, in the following sense: there are sofic approximations $(\sigma_n')$ and $(\tau_n')$ of $\Gamma$ such that $(\sigma_n)$ is equivalent to $(\sigma_n')$, $(\sigma_n')$ is branched covering of $(\tau_n')$ and $(\tau_n')$ is equivalent to $(\tau_n)$.

\begin{remark} \label{rem2} To give a better understanding of the notion of being essentially a branched covering, we record two equivalent formulations of it. 
Let $\sigma_n:\Gamma\rightarrow\text{Sym}(X_n)$ and $\tau_n:\Gamma\rightarrow\text{Sym}(Y_n)$ be asymptotic homomorphisms of a countable group $\Gamma$. Then the following conditions are equivalent:
\begin{enumerate}[label=(\roman*)]
\item $(\sigma_n)$ is essentially a branched covering of $(\tau_n)$.
\item there are
onto,  $d_n$-to-one maps  $\theta_n:X_n\rightarrow Y_n$, for some $d_n\in\mathbb N$, such that we have $\lim\limits_{n\rightarrow\infty}\text{d}_{\text{H}}(\theta_n\circ\sigma_n(g),\tau_n(g)\circ\theta_n)=0$, for every $g\in\Gamma$.
\item $(\sigma_n)$ is equivalent to a branched covering $(\sigma_n')$ of $(\tau_n)$ (i.e., one can take $\tau_n'=\tau_n$ in (i)).
\end{enumerate}
It is clear that (i) $\Rightarrow$ (ii) and (iii) $\Rightarrow$ (i). 
That (ii) $\Rightarrow$ (iii) is a consequence of the following fact: if $X, Y$ are finite sets, $\sigma\in\text{Sym}(X)$, $\tau\in\text{Sym}(Y)$, and $\theta:X\rightarrow Y$ an onto, $d$-to-one map, for $d\in\mathbb N$, then there is $\sigma'\in\text{Sym}(X)$ such that $\theta\circ\sigma'=\tau\circ\theta$ and $\text{d}_{\text{H}}(\sigma',\sigma)\leq\text{d}_{\text{H}}(\theta\circ\sigma,\tau\circ\theta)$. 
\end{remark}

Our main  result settles Bowen's question in the negative. More precisely, we prove the following:

\begin{theorem}\label{main}
Let $\Gamma=\mathbb F_m\times\mathbb F_k$, for some integers $m,k\geq 2$.
Then $\Gamma$ admits a sofic approximation $\sigma_n:\Gamma\rightarrow \emph{Sym}(X_n)$  with the following property: there are no homomorphisms $\tau_n:\Gamma\rightarrow\emph{Sym}(Y_n)$ and maps $\theta_n:X_n\rightarrow Y_n$, for some finite sets $Y_n$, such that 
\begin{enumerate}[label=\emph{(\alph*)}]
\item $\lim\limits_{n\rightarrow\infty}\emph{d}_{\emph{H}}(\theta_n\circ\sigma_n(g),\tau_n(g)\circ\theta_n)=0$, for all $g\in\Gamma$, and
\item $\lim\limits_{n\rightarrow\infty}\emph{d}_{\emph{H}}(\theta_n\circ\sigma_n(g),\theta_n)=1$, for all $g\in\Gamma\setminus\{e\}$.
\end{enumerate}
\end{theorem}

Theorem \ref{main} implies that $(\sigma_n)$ is not essentially a branched covering of a sofic approximation by homomorphisms. This follows by using 
Remark \ref{rem2} ((i) $\Rightarrow$ (ii)) and noting that if each $\theta_n:X_n\rightarrow Y_n$ is $d_n$-to-one, for some $d_n\in\mathbb N$, and $(\tau_n)$ is a sofic approximation  of $\Gamma$, then (a) implies (b).

 Theorem \ref{main} strengthens part of \cite[Theorem D]{Io19b}.  More precisely, \cite[Theorem D]{Io19b} shows that $\Gamma=\mathbb F_m\times\mathbb F_k$ is not weakly very flexibly stable, for any  integers $m,k\geq 2$, in the sense of \cite[Definition 1.6]{Io19b}. This amounts to the existence of a sofic approximation $\sigma_n:\Gamma\rightarrow\text{Sym}(X_n)$ with the following property: $(\star)$ there are no finite sets $Y_n$, homomorphisms $\tau_n:\Gamma\rightarrow\text{Sym}(Y_n)$ and one-to-one maps $\theta_n:X_n\rightarrow Y_n$ such that $\lim\limits_{n\rightarrow\infty}\text{d}_{\text{H}}(\theta_n\circ\sigma_n(g),\tau_n(g)\circ\theta_n)=0$, for every $g\in\Gamma$.

 As we explain in the comments below, the sofic approximation $(\sigma_n)$ of $\Gamma$  from the hypothesis of Theorem \ref{main} is constructed following the strategy introduced in \cite{Io19b}. As such, results from \cite{Io19b} readily imply that  $(\sigma_n)$ satisfies $(\star)$.
 The main novelty in the proof of Theorem \ref{main} consists of showing that any maps $\theta_n:X_n\rightarrow Y_n$ as in its statement must be  ``asymptotically one-to-one".
Moreover, we prove that if $\tau_n:\Gamma\rightarrow\text{Sym}(Y_n)$ are arbitrary maps, then any maps $\theta_n:X_n\rightarrow Y_n$ which satisfy conditions (a) and (b) from Theorem \ref{main} must be asymptotically one-to-one. This implies that $(\sigma_n)$ is a minimal sofic approximation of $\Gamma$, in the sense that it is not equivalent to a proper  (i.e., one satisfying $d_n>1$, for every $n$) branched covering of any sofic approximation of $\Gamma$.

 \subsection*{Comments on the proof of Theorem \ref{main}}
 We end the introduction with an informal outline of the proof of our main result. If a group $\Gamma$ satisfies the conclusion of Theorem \ref{main} then any group containing it as a finite index subgroup also does (see Lemma \ref{fin}). Therefore, it suffices to prove Theorem \ref{main} 
 when $m\geq 5$ and $k\geq 3$. Fix a free decomposition $\mathbb F_m=\mathbb F_{m-1}*\mathbb Z$. Thus, we have $\Gamma=\mathbb F_m\times\mathbb F_k=(\mathbb F_{m-1}*\mathbb Z)\times\mathbb F_k$. The proof of Theorem \ref{main} is divided between Sections \ref{3} and \ref{4}:
 \begin{enumerate}
 
\item In  Section \ref{3}, we use the work \cite{ALW01}  of Alon, Lubotzky and Widgerson  who proved that expansion is not a group property.
This allows us to define a sequence of finite groups $G_p$ (indexed over primes $p\equiv 1\pmod{3}$) together with onto homomorphisms $\varphi_p:\mathbb F_{m-1}\rightarrow G_p$ and $\rho_p:\mathbb F_k\rightarrow G_p$ such that  $\mathbb F_{m-1}$ has property $(\tau)$ with respect to $\{\ker(\varphi_p)\}$, while $\mathbb F_k$ does not have property ($\tau)$ with respect to $\{\ker(\rho_p)\}$. (A key property of $G_p$ is that it has only one non-trivial normal subgroup. To the best of our knowledge, it is unknown if one can find such groups $G_p$ which are simple; if this were the case, then the proof could be simplified considerably.)
Following closely \cite{Io19b}  we then construct an asymptotic homomorphism $\sigma_p:\mathbb F_m\times\mathbb F_k\rightarrow\text{Sym}(G_p)$  which satisfies $(\star)$ and that $\sigma_p(g,h)x=\varphi_p(g)x\rho_p(h)^{-1}$, for every $g\in\mathbb F_{m-1}, h\in\mathbb F_k$ and $x\in G_p$. Note, however, that the asymptotic homomorphism $(\sigma_p)$ is not a sofic approximation. 
 
 \item We begin Section \ref{4} by augmenting the construction of $(\sigma_p)$ to get a sofic approximation $\widetilde\sigma_p:\mathbb F_m\times\mathbb F_k\rightarrow\text{Sym}(\widetilde G_p)$ which inherits the properties of $(\sigma_p)$ listed above. 
The rest of Section \ref{4} is devoted to proving that $(\widetilde\sigma_p)$ verifies the conclusion of Theorem \ref{main}. 
Assume by contradiction that there are homomorphisms $\tau_p:\mathbb F_m\times\mathbb F_k\rightarrow\text{Sym}(Y_p)$, for some finite sets $Y_p$, and maps $\theta_p:\widetilde{G_p}\rightarrow Y_p$ which satisfy conditions  (a) and (b) from Theorem \ref{main}. Condition (b) implies that the partition $\{\theta_p^{-1}(\{y\})\mid y\in Y_p\}$ of $\widetilde G_p$ is $\widetilde\sigma_p(\mathbb F_m\times\mathbb F_k)$-asymptotically invariant. By combining the property $(\tau)$ assumption with a result from \cite{Io19a} (see Section \ref{2}), we deduce that the partition $\{\theta_p^{-1}(\{y\})\mid y\in Y_p\}$ is asymptotically equal to the coset partition $\{gN_p\mid g\in \widetilde G_p\}$ of $\widetilde G_p$, for some normal subgroup $N_p\trianglelefteq\widetilde G_p$. Some additional work, which uses condition (b), allows us to conclude that $N_p=\{e\}$, and thus $\theta_p$ is asymptotically one-to-one. This however contradicts the fact that $(\widetilde\sigma_p)$ satisfies $(\star)$. 
 \end{enumerate}

\subsection*{Acknowledgements.} I am grateful to Lewis Bowen for helpful discussions clarifying his question answered here.

\section{Property $(\tau)$ and almost invariant partitions}\label{2}
 This section is devoted to a technical lemma which will be needed in the proof of our main theorem. 
Let $\Gamma$ be a finitely generated group, $S$ be a finite set of generators of $\Gamma$ and $\{\Gamma_n\}_{n=1}^{\infty}$ be a sequence of finite index  normal subgroups. Denote $G_n=\Gamma/\Gamma_n$ and let $p_n:\Gamma\rightarrow G_n$ be the quotient homomorphism. The following lemma asserts that if $\Gamma$ has property $(\tau)$ with respect to $\{\Gamma_n\}_{n=1}^{\infty}$, then any partition of $G_n$ which is almost invariant under the left multiplication action of $\Gamma$ must essentially come from the left cosets of a subgroup of $G_n$.

Recall that $\Gamma$ is said to have {\it property ($\tau$)} with respect to $\{\Gamma_n\}_{n=1}^{\infty}$ if $\inf_n\kappa(G_n,p_n(S))>0$ \cite{Lu94}.
Here, given a finite group $G$ and a set of generators $T\subset G$, the {\it Kazhdan constant} $\kappa(G,T)$ denotes the largest constant
 $\kappa>0$ such that $\kappa\cdot \|\xi\|\leq\max_{g\in T}\|\pi(g)\xi-\xi\|$, for every $\xi\in\mathcal H$ and unitary representation $\pi:G\rightarrow\text{U}(\mathcal H)$ of $G$ on a Hilbert space $\mathcal H$ which has no non-zero invariant vectors. We record the following remark which will be needed in the proof of Lemma \ref{ALW01}.

\begin{remark}\label{arb}
Let $\pi:G\rightarrow\text{U}(\mathcal H)$ be a unitary representation of a finite group $G$. Let $P$ be the orthogonal projection from $\mathcal H$ onto  the closed subspace  $\mathcal H^G$ of $\pi(G)$-invariant vectors and $\xi\in\mathcal H$. Then $\max_{g\in G}\|\pi(g)\xi-\xi\|\leq 2\cdot \|\xi-P(\xi)\|$. Since the restriction of $\pi$ to $\mathcal H\ominus\mathcal H^G$ has no non-zero invariant vectors, we get that $\kappa(G,T)\cdot\|\xi-P(\xi)\|\leq\max_{g\in T}\|\pi(g)\xi-\xi\|$ and further that\begin{equation}\label{arbrep}\text{$\frac{\kappa(G,T)}{2}\cdot\max_{g\in G}\|\pi(g)\xi-\xi\|\leq\max_{g\in T}\|\pi(g)\xi-\xi\|$, for every $\xi\in\mathcal H$.}\end{equation} \end{remark}

\begin{lemma}\emph{\cite{Io19a}} 
\label{Io18}
In the above setting, assume that $\Gamma$ has property $(\tau)$ with respect to $\{\Gamma_n\}_{n=1}^{\infty}$.
For every $n$, let $\{X_{n,k}\}_{k=1}^{d_n}$ be a partition of $G_n$, for some $d_n\in\mathbb N$. Assume that for every $n$ and $g\in\Gamma$, there exists a permutation $\sigma_{n,g}$ of $\{1,...,d_n\}$ such that $\lim\limits_{n\rightarrow\infty}\frac{1}{|G_n|}\cdot\sum_{k=1}^{d_n}|gX_{n,k}\triangle X_{n,\sigma_{n,g}(k)}|=0.$

Then for every $n$ we can find a subgroup $H_n<G_n$, a set $S_n\subset\{1,...,d_n\}$ and a one-to-one map $\omega_n:S_n\rightarrow G_n/H_n$ such that $$\text{$\lim\limits_{n\rightarrow\infty}\frac{1}{|G_n|}\cdot\sum_{k\in S_n}|X_{n,k}\triangle \omega_n(k)H_n|=0$\;\;\; and \;\;\; $\lim\limits_{n\rightarrow\infty}\frac{1}{|G_n|}\sum_{k\notin S_n}|X_{n,k}|=0$}.$$
\end{lemma}
 This result is a consequence of the proof of \cite[Theorem A]{Io19a}. For the reader's convenience, we indicate briefly how  the proof of \cite[Theorem A]{Io19a} can be adapted to prove Lemma \ref{Io18}.

{\it Proof.}
For every $n$, let $\pi_n:\Gamma\rightarrow \text{U}(\ell^2(G_n\times G_n))$ be the unitary representation associated to the action $\Gamma\curvearrowright G_n\times G_n$ given by $g\cdot (x,y)=(gx,gy)$, and define the unit vector $$\eta_n=\sum_{k=1}^{d_n}\frac{1}{\sqrt{|X_{n,k}}|} 
\;{\bf 1}_{X_{n,k}\times X_{n,k}}\in\ell^2(G_n\times G_n).$$

We claim that $\|\pi_n(g)\eta_n-\eta_n\|_2\rightarrow 0$, for every $g\in\Gamma$. To this end, fix $g\in\Gamma$. Then the hypothesis implies that $\frac{1}{|G_n|}\cdot\sum_{k=1}^{d_n}|gX_{n,k}\cap X_{n,\sigma_{n,g}(k)}|\rightarrow 1$ and thus $\frac{1}{|G_n|}\cdot\sum_{k=1}^{d_n}\sqrt{|X_{n,k}|\cdot |X_{n,\sigma_{n,g}(k)}|}\rightarrow 1$. Using a direct computation and the Cauchy-Schwarz inequality we derive that \begin{align*}
\langle\pi_n(g)\eta_n,\eta_n\rangle&=\sum_{k,l=1}^{d_n}\frac{1}{\sqrt{|X_{n,k}|\cdot |X_{n,l}|}}\;|gX_{n,k}\cap X_{n,l}|^2\\ &\geq \sum_{k=1}^{d_n}\frac{1}{\sqrt{|X_{n,k}|\cdot |X_{n,\sigma_{n,g}(k)}|}}\; |gX_{n,k}\cap X_{n,\sigma_{n,g}(k)}|^2\\&\geq\frac{\sum_{k=1}^{d_n}|gX_{n,k}\cap X_{n,\sigma_{n,g}(k)}|}{\sum_{k=1}^{d_n}\sqrt{|X_{n,k}|\cdot |X_{n,\sigma_{n,g}(k)}|}}.\end{align*} Thus, $\liminf_{n\rightarrow\infty}\langle\pi_n(g)\eta_n,\eta_n\rangle\geq 1$ and since $\|\eta_n\|_2=1$, we conclude that  $\|\pi_n(g)\eta_n-\eta_n\|_2\rightarrow 0$. 

Since $\Gamma$ has property $(\tau)$ with respect to $\{\Gamma_n\}_{n=1}^{\infty}$, we get that $\kappa:=\inf_n\kappa(G_n,p_n(S))>0$.
By \cite[Lemma 2.5]{Io19b} we deduce that $\sup_{g\in\Gamma}\|\pi_n(g)\eta_n-\eta_n\|_2\leq (2/\kappa)\cdot\max_{g\in S}\|\pi_n(g)\eta_n-\eta_n\|_2$, for every $n$. In combination with the above  it follows that $\sup_{g\in\Gamma}\|\pi_n(g)\eta_n-\eta_n\|_2\rightarrow 0$. Thus, we can find positive real numbers $\delta_n$ such that $\delta_n\rightarrow 0$ and 
$\sup_{g\in\Gamma}\|\pi_n(g)\eta_n-\eta_n\|_2^2<2\delta_n$, for every $n$.

Let $n$ large enough such that $\delta_n<10^{-12}$. Then \begin{equation}\label{prof}\text{$\sum_{k,l=1}^{d_n}\frac{1}{\sqrt{|X_{n,k}|\cdot |X_{n,l}|}}\;|gX_{n,k}\cap X_{n,l}|^2=\langle\pi_n(g)\eta_n,\eta_n\rangle>1-\delta_n$, for every $g\in\Gamma$}.\end{equation}
By using \eqref{prof} and applying verbatim the second part of the proof of \cite[Theorem A]{Io19a}, we can find a subgroup $H_n<G_n$, a nonempty subset $S_n\subset\{1,...,d_n\}$ and a map $\omega_n:S_n\rightarrow G_n/H_n$ such that $\sum_{k\in S_n}|X_{n,k}|\geq (1-\sqrt{\delta_n})\cdot |G_n|$, and 
$|X_{n,k}\triangle\omega_n(k)H_n|\leq 506\sqrt[4]{\delta_n}\cdot |X_{n,k}|$, for every $k\in S_n$. 

If $k,l\in S_n$ and $k\not=l$, then using that $X_{n,k}\cap X_{n,l}=\emptyset$, we get that \begin{align*}|\omega_n(k)H_n\triangle \omega_n(l)H_n|&\geq |X_{n,k}\triangle X_{l,n}|-|X_{n,k}\triangle \omega_n(k)H_n|-|X_{n,l}\triangle \omega_n(l)H_n|\\& \geq (1- 506\sqrt[4]{\delta_n})\cdot (|X_{n,k}|+|X_{n,l}|).\end{align*}
Since $ 506\sqrt[4]{\delta_n}<1$, we derive that $\omega_n(k)H_n\triangle \omega_n(l)H_n\not=\emptyset$.
This implies that the map $\omega_n$ is one-to-one and the conclusion follows.
\hfill$\blacksquare$
\section{Construction of asymptotic homomorphisms}\label{3} 
In this section, we establish two ingredients that will be needed in the proof of our main theorem. To explain this, fix integers $m\geq 5$ and $k\geq 3$, and denote $\Gamma=\mathbb F_{m-1}$ and $\Lambda=\mathbb F_k$. In the first part of this section, given a prime $p$ with $p\equiv 1 \pmod{3}$, we construct a finite group $G_p$ and homomorphisms $\varphi_p:\Gamma\rightarrow G_p$, $\rho_p:\Lambda\rightarrow G_p$ with various special properties.  In the second part of this section, we follow closely \cite[Section 6]{Io19b} to construct an asymptotic homomorphism 
$\sigma_p:(\Gamma*\mathbb Z)\times\Lambda\rightarrow \text{Sym}(G_p)$ such that $\sigma_p(g,h)x=\varphi_p(g)x\rho_p(h)^{-1}$, for all $g\in\Gamma,h\in\Lambda$ and $x\in G_p$. 

\subsection{A group theoretic construction} 
In \cite{ALW01}, Alon, Lubotzky and Widgerson showed that expansion is not a group property. Thus, they introduced a method of constructing sequences of finite groups $\{G_n\}_{n=1}^{\infty}$ and generating sets $S_n,T_n$ of fixed cardinality ($|S_n|=m$, $|T_n|=k$) such that  the Cayley graphs of $G_n$ are expanders with respect to $S_n$ but not with respect to $T_n$.
Equivalently, there are onto homomorphisms $p_n:\mathbb F_m\rightarrow G_n$, $q_n:\mathbb F_{k}\rightarrow G_n$ such that $\mathbb F_m$ has property $(\tau)$ with respect to $\{\ker(p_n)\}_{n=1}^{\infty}$ but $\mathbb F_{k}$ does not have property $(\tau)$ with respect to $\{\ker(q_n)\}_{n=1}^{\infty}$.

The proof of our main theorem relies on a particular case of the construction of \cite{ALW01}.
Let $p$ be a prime with $p\equiv 1 \pmod{3}$. Denote by $\text{P}^1(F_p)=F_p\cup\{\infty\}$ the projective line over the field $F_p$ with $p$ elements.
Consider the action of $\text{PSL}_2(F_p)=\text{SL}_2(F_p)/\{\pm\text{I}\}$ on $\text{P}^1(F_p)$ by  linear fractional transformations: $$\begin{pmatrix}a&b\\c&d\end{pmatrix}\cdot x=\frac{ax+b}{cx+d}.$$ 

Further, we consider the vector space $F_3^{\text{P}^1(F_p)}$ over $F_3$, and the permutational representation of $\text{PSL}_2(F_p)$ on $F_3^{\text{P}^1(F_p)}$ given by $g\cdot x=(x_{g^{-1}\cdot i})_{i\in\text{P}^1(F_p)}$, for every $g\in\text{PSL}_2(F_p)$ and  $x=(x_i)_{i\in \text{P}^1(F_p)}$. We identify $F_3^{\text{P}^1(F_p)}$ with $F_3^{p+1}$ using a fixed bijection $\text{P}^1(F_p)\mapsto \{1,...,p+1\}$ which sends $\infty$ to $p+1$. 
We continue by introducing the following:

\begin{notation}\label{G_p}
We denote $A_p=\{(x_i)\in F_3^{p+1}\mid\sum_{i=1}^{p+1}x_i=0\}$ and $H_p=\text{PSL}_2(F_p)$. Then $A_p$ is an $H_p$-invariant subspace of $F_3^{p+1}$ with $|A_p|=|F_3^{p+1}|/3=3^p$. We denote $G_p=A_p\rtimes H_p.$
\end{notation}

We next record the following elementary result, whose proof we include for completeness.

\begin{lemma}\label{irred}
Let $K$ be a group and $N<G_p\times K$ be a subgroup which is normalized by $G_p\times\{e\}$. Then $N$ is equal to $\{e\}\times L$, $A_p\times L$ or $G_p\times L$, for some subgroup $L<K$.
\end{lemma}

{\it Proof.}
First, we claim that if $N<A_p$ is a subgroup which is normalized by $H_p$, then $N=\{e\}$ or $N=A_p$.  To this end, suppose that $N$ contains an element $x=(x_i)_{i=1}^{p+1}$ not equal to $(0,..,0)$. Since $N$ is $H_p$-invariant and $H_p$ acts transitively on $\text{P}^1(F_p)$, we may assume that $y:=x_{p+1}\not=0$.
 Since the subgroup $U_p=\{\begin{pmatrix}1&a\\0&1 \end{pmatrix}\mid a\in F_p\}$ of $H_p$ fixes $\infty\in \text{P}^1(F_p)$ and acts transitively on $F_p$, we get
   $$\sum_{g\in U_p}g\cdot x=(\sum_{i=1}^px_i,...,\sum_{i=1}^px_i,px_{p+1})=(-y,...,-y,py)\in N.$$
 Since $y\not=0$, we get that $(-1,...,-1,p)\in N$. 
  Let $e_k=(x_ {k,i})_{i=1}^{p+1}$, where $x_{k,i}=-1$ for $i\not=k$ and $x_{k,k}=p$.
 Since $H_p$ acts transitively on $\text{P}^1(F_p)$, we get that $e_k\in N$, for all $1\leq k\leq p$. Since $3\nmid p+1$, the vectors $(e_k)_{k=1}^p\in A_p$ are linearly independent over $F_3$. As $\dim A_p=p$, we get that $N=A_p$.
 
 Second, we claim that if $N<G_p$ is a normal subgroup, then $N=\{e\}$, $N=A_p$ or $N=G_p$. Let $\rho:G_p\rightarrow H_p$ be the quotient homomorphism. Then $\rho(N)<H_p$ is a normal subgroup and since $H_p$ is a simple group,  $\rho(N)=\{e\}$ or $\rho(N)=H_p$. If $\rho(N)=\{e\}$, then $N<A_p$ and the first claim implies that $N=\{e\}$ or $N=A_p$. It remains to analyze the case when $\rho(N)=H_p$.
We first show that $N\cap A_p\not=\{e\}$. Assume by contradiction that $N\cap A_p=\{e\}$ and let $a\in A_p\setminus\{e\}$. Since $\rho(N)=H_p$, for any $h\in H_p$, there is $b\in A_p$ such that $bh\in N$. Since $N<G_p$ is normal, $abha^{-1}\in N$ and thus $aha^{-1}h^{-1}=(abha^{-1})(bh)^{-1}\in N\cap A_p$. Hence $aha^{-1}h^{-1}=e$, for every $h\in H_p$, which contradicts that $a\not=e$. Finally, if $N\cap A_p\not=\{e\}$, then since $N\cap A_p<A_p$ is normalized by $H_p$, the first claim implies that $N\cap A_p=A_p$ and hence $N\supset A_p$. Since $\rho(N)=H_p$, it follows that $N=G_p$. 
 
Let $N<G_p\times K$ be a subgroup which is normalized by $G_p\times\{e\}$. The second claim implies that $N\cap (G_p\times\{e\})$ is equal to $\{e\}$, $A_p\times\{e\}$ or $G_p\times\{e\}$. Note that if $(g,k)\in N$, for some $g\in G_p$ and $k\in K$, then $(ghg^{-1}h^{-1},e)=(g,k)(h,e)(g,k)^{-1}(h,e)^{-1}\in N\cap (G_p\times\{e\})$, for every $h\in G_p$. If $N\cap (G_p\times\{e\})=\{e\}$, it follows that $N\subset \{e\}\times K$, thus $N=\{e\}\times L$, for some subgroup $L<K$. If $N\cap (G_p\times\{e\})=A_p\times\{e\}$, we get that if $(g,k)\in N$, then $ghg^{-1}h^{-1}\in A_p$, for every $h\in G_p$, and thus $g\in A_p$. Hence $A_p\times\{e\}\subset N\subset A_p\times K$, which implies that $N=A_p\times L$, for some subgroup $L<K$. Finally, if $N\cap (G_p\times\{e\})=G_p\times\{e\}$, then $N=G_p\times L$, for some subgroup $L<K$.
\hfill$\blacksquare$

In addition to the notation from \ref{G_p}, throughout the rest of this paper we will use the following:
\begin{notation} \label{setting}
Given a prime $p$ with $p\equiv 1\pmod {3}$, we fix a prime $r_p>p$, denote $K_p=\text{PSL}_2(F_{r_p})$ and let $\psi_p:\text{PSL}_2(\mathbb Z)\rightarrow K_p$  be the quotient homomorphism. We denote $\widetilde G_p=G_p\times K_p$.
\end{notation}

The following result combines \cite{ALW01} with a spectral result gap result from \cite{BV12}.

\begin{lemma}\label{ALW01} Let $\Gamma=\mathbb F_{m-1}$, for  $m\geq 5$. View $\Gamma$ as a subgroup of $\emph{PSL}_2(\mathbb Z)$. 
Then for any large enough prime $p$ with $p\equiv 1\pmod {3}$, there is an onto homomorphism $\varphi_p:\Gamma\rightarrow G_p$ such that $\widetilde\varphi_p:\Gamma\rightarrow\widetilde G_p$ given by $\widetilde\varphi_p(g)=(\varphi_p(g),\psi_p(g))$ is onto and $\Gamma$ has property $(\tau)$ with respect to $\{\ker(\widetilde\varphi_p)\}_p$.
\end{lemma}

{\it Proof.}  Let   $a_1,...,a_{m-1}$  be free generators of $\Gamma$ and $p$ be a prime with $p\equiv 1 \pmod{3}$.
Denote by $\xi_p:\text{PSL}_2(\mathbb Z)\rightarrow H_p$ the quotient homomorphism and let $\eta_p:\text{PSL}_2(\mathbb Z)\rightarrow H_p\times K_p$ be the homomorphism given by $\eta_p(g)=(\xi_p(g),\psi_p(g))$.
Since $H_p\times K_p=\text{PSL}_2(F_p)\times\text{PSL}_2(F_{r_p})\cong\text{PSL}_2(\mathbb Z/pr_p\mathbb Z)$ and $\langle a_1,...,a_{m-3}\rangle\cong\mathbb F_{m-3}$ is a non-amenable subgroup of $\text{PSL}_2(\mathbb Z)$ (as $m-3\geq 2$), it follows that for large enough $p$ we have $\eta_p(\langle a_1,...,a_{m-3}\rangle)=H_p\times K_p$.

By applying \cite[Theorem 1]{BV12}, we conclude that \begin{equation}\label{spectral2} \kappa_1:=\inf_p\kappa(H_p\times K_p,\{\eta_p(a_1),....,\eta_p(a_{m-3})\})>0. \end{equation}
For $w\in A_p$, we denote by $w^{H_p}=\{h(w)=hwh^{-1}\mid h\in H_p\}$ the orbit of $w$ under the action  of $H_p$. 
By Lemma \ref{irred}, the permutational representation of $H_p$ on $A_p\subset F_3^{p+1}$ is irreducible. Thus, by applying \cite[Theorem 3.1]{ALW01}, we can find $v_1(p),v_2(p)\in A_p\setminus\{e\}$ such that
\begin{equation}\label{spectral}
\kappa_2:=\inf_p\kappa(A_p, v_1(p)^{H_p}\cup v_2(p)^{H_p})>0.\end{equation}
Define a homomorphism $\varphi_p:\Gamma\rightarrow G_p=A_p\rtimes H_p$ by letting $\varphi_p(a_i)=\xi_p(a_i)$, for $1\leq i\leq m-3$, $\varphi_p(a_{m-2})=v_1(p)\xi_p(a_{m-2})$ and $\varphi_p(a_{m-1})=v_2(p)\xi_p(a_{m-1})$. Then $\widetilde\varphi_p:\Gamma\rightarrow\widetilde G_p$ given by $\widetilde\varphi_p(g)=(\varphi_p(g),\psi_p(g))$ is onto. Indeed, since $\widetilde\varphi_p(a_i)=\eta_p(a_i)$, for all $1\leq i\leq m-3$, we get that
$\widetilde\varphi_p(\Gamma)$ contains $\eta_p(\langle a_1,...,a_{m-3}\rangle)=H_p\times K_p$. Thus, $\widetilde\varphi_p(\Gamma)$ also contains $(v_1(p),e)\in (A_p\setminus\{e\})\times\{e\}$. Since $A_p$ has no proper non-trivial $H_p$-invariant subgroup by Lemma \ref{irred}, we derive that $\widetilde\varphi_p$ is onto. In particular, $\varphi_p$ is onto.

Moreover, combining \eqref{spectral2} and \eqref{spectral} as in \cite{ALW01} implies that $\inf_p\kappa(\widetilde G_p,\{\widetilde\varphi_p(a_1),...,\widetilde\varphi_p(a_m)\})>0$.  
To justify this, put $\kappa:=\frac{\min\{\kappa_1,\kappa_2\}}{2}>0$. Let $\pi:\widetilde G_p\rightarrow\text{U}(\mathcal H)$ be a unitary representation with no non-zero invariant vectors and $\xi\in\mathcal H$.
For $F\subset\widetilde G_p$, let $\Delta(F)=\max_{g\in F}\|\pi(g)\xi-\xi\|$ and note that $\Delta(F)\leq \Delta(F_1)+\Delta(F_2)$, whenever $F\subset F_1 F_2$.
 By combining \eqref{arbrep} with \eqref{spectral2} and \eqref{spectral} we get that
\begin{equation}\label{una}\text{$\kappa\cdot \Delta(H_p\times K_p)\leq \Delta(\{\varphi_p(a_i)\mid 1\leq i\leq m-3\})$\;\; and}
\end{equation}
\begin{equation}\label{doua}
\kappa\cdot\Delta(A_p\times\{e\})\leq\Delta((v_1(p)^{H_p}\cup v_2(p)^{H_p})\times\{e\}).
\end{equation} 
Since $(v_1(p),e)\in\widetilde\varphi_p(a_{m-2})(H_p\times K_p)$  we derive that $v_1(p)^{H_p}\times\{e\}\subset  (H_p\times K_p)\widetilde\varphi_p(a_{m-2})(H_p\times K_p)$.
Thus, we have that \begin{equation}\label{trei} \Delta(v_1(p)^{H_p}\times\{e\})\leq \Delta(\{\widetilde\varphi_p(a_{m-2})\})+2\cdot\Delta(H_p\times K_p).\end{equation}Similarly, we get that \begin{equation}\label{patru} \Delta(v_2(p)^{H_p}\times\{e\})\leq \Delta(\{\widetilde\varphi_p(a_{m-1})\})+2\cdot\Delta(H_p\times K_p).\end{equation} Since $\widetilde G_p=(A_p\times\{e\})(H_p\times K_p)$, by combining \eqref{una}, \eqref{doua}, \eqref{trei} and \eqref{patru} it is easy to see that \begin{equation}\label{cinci}\frac{\kappa^2}{2(\kappa+1)}\cdot\Delta(\widetilde G_p)\leq\Delta(\{\widetilde\varphi_p(a_i)\mid 1\leq i\leq m-1\}).\end{equation} Since $\pi$ has no non-zero invariant vectors we have that $\|\xi\|\leq\Delta(\widetilde G_p)$. Together with \eqref{cinci}, this implies that $\inf_p\kappa(\widetilde G_p,\{\widetilde\varphi_p(a_1),...,\widetilde\varphi_p(a_{m-1})\})\geq\frac{\kappa^2}{2(\kappa+1)}>0$. Hence, $\Gamma$ has property $(\tau)$ with respect to $\{\ker(\widetilde\varphi_p)\}_p$, which finishes the proof of the lemma.
\hfill$\blacksquare$

\begin{lemma}\label{ALW01II}
Let $b_1,...,b_k$ be free generators of $\Lambda=\mathbb F_k$, for $k\geq 3$. Then there is $C>0$ such that for every large enough prime $p$ with $p\equiv 1 \pmod{3}$ there is an onto homomorphism $\rho_p:\Lambda\rightarrow G_p$ and a set $T_p\subset G_p$ satisfying that  $\frac{1}{243}\leq\frac{|T_p|}{|G_p|}\leq\frac{1}{3}$ and $|T_p\rho_p(b_j)\triangle T_p|\leq\frac{C}{\sqrt{p}}|G_p|$, for every $1\leq j\leq k$. Moreover, if $h\in\Lambda\setminus\{e\}$, then $\rho_p(h)\not=e$, for every large enough $p$.
\end{lemma}

{\it Proof.} Let $p$ be a prime with $p\equiv 1\pmod{3}$ and put $v(p)=(1,-1,0,...,0)\in A_p$.
For $x=(x_i)\in A_p$ and $j\in\{0,1,2\}$, denote $n_j(x)=|\{i\mid x_i=j\pmod{3}\}|$. We define \begin{equation}\label{S_p}S_p=\{x\in A_p\mid\text{$n_1(x)>n_0(x)+2$ and $n_1(x)>n_2(x)+2$}\}.\end{equation}
First, we claim that \begin{equation}\label{bounds}\frac{1}{243}\leq\frac{|S_p|}{|A_p|}\leq\frac{1}{3}.
\end{equation}Since $p\equiv 1\pmod{3}$, we  have $a_1=(2,1,...,1)\in A_p$ and $a_2=(1,2,...,2)\in A_p$. Since the sets $S_p, a_1+S_p, a_2+S_p$ are pairwise disjoint, it follows that $|S_p|\leq\frac{|A_p|}{3}$.
On the other hand, if we let $R_p=\{(x_i)\in F_3^{p-4}\mid \text{$n_1(x)\geq n_0(x)$ and $n_1(x)\geq n_2(x)$}\}$, then for every $x=(x_i)\in R_p$, there is $\tilde x=(\tilde x_i)\in S_p$ with $\tilde x_i=x_i$, for all $1\leq i\leq p-4$. Thus, $|S_p|\geq |R_p|\geq \frac{|F_3^{p-4}|}{3}=\frac{|A_p|}{243}$, proving \eqref{bounds}.

Second, we claim that there is a constant $C>0$ such that
\begin{equation}\label{eqC}\text{$|(v(p)+S_p)\triangle S_p|\leq \frac{C}{\sqrt{p}}|A_p|$,\;\;\; for every $p$.}
\end{equation}
To this end, note that Stirling's formula implies that there is a constant $c>0$ such that
\begin{equation}\label{stirling}
{n \choose k}\leq {n\choose\lfloor\frac{n}{2}\rfloor}\leq c\cdot\frac{2^n}{\sqrt{n}}, \;\;\text{for every $n\geq k\geq 0$}.
\end{equation}
Since $S_p\setminus (v(p)+S_p)\subset\{x\in F_p^{p+1}\mid\text{$n_0(x)+5\geq n_1(x)>n_0(x)+2$ or $n_2(x)+5\geq n_1(x)> n_2(x)+2$}\}$ and $n_0(x),n_2(x)<\frac{p+1}{2}$, for every $x\in S_p$, 
 by using \eqref{stirling} we get that \begin{align*}|S_p\setminus (v(p)+S_p)|&\leq 2\cdot\sum_{\substack{{n_0+n_1+n_2=p+1}\\{0\leq n_0<\frac{p+1}{2}}\\{ n_2+2<n_1\leq n_2+5}}}\frac{(p+1)!}{n_0!n_1!n_2!}\\&=2\cdot\sum_{0\leq n_0<\frac{p+1}{2}}{p+1 \choose n_0}\sum_{\substack{{n_1+n_2=p+1-n_0}\\ {n_2+2<n_1\leq n_2+5}}}{p+1-n_0\choose n_1}\\&\leq 6c\cdot\sum_{0\leq n_0<\frac{p+1}{2}}{p+1 \choose n_0}\frac{2^{p+1-n_0}}{\sqrt{p+1-n_0}}\\&\leq 6c\cdot\sqrt{\frac{2}{p+1}}\cdot\sum_{0\leq n_0<\frac{p+1}{2}}{p+1 \choose n_0}2^{p+1-n_0}\\&\leq 6c\cdot\sqrt{\frac{2}{p+1}}\cdot\frac{3^{p+1}}{2}
\end{align*}
Since $|A_p|=3^p$, this proves \eqref{eqC}.

Let $\xi_p:\text{PSL}_2(\mathbb Z)\rightarrow H_p=\text{PSL}_2(F_p)$ be the quotient homomorphism. 
 Since $\langle b_1,...,b_{k-1}\rangle\cong \mathbb F_{k-1}$ is a non-amenable subgroup of $\text{PSL}_2(\mathbb Z)$, for every large enough $p$ we have $\xi_p(\langle b_1,...,b_{k-1}\rangle)=H_p$.  
We define a homomorphism $\rho_p:\Lambda\rightarrow G_p$  by  $\rho_p(b_j)=\xi_p(b_j)$, for every $1\leq j\leq k-1$, and $\rho_p(b_k)=\xi_p(b_k)v(p)$. Since $A_p$ has no proper non-trivial $H_p$-invariant subgroup by Lemma \ref{irred}, it follows that $\rho_p$ is onto. 

Next, note that $S_p$ is $H_p$-invariant and let $T_p=H_p\cdot S_p=S_p\cdot H_p\subset G_p$. Then \eqref{bounds} and \eqref{eqC} imply that $\frac{|G_p|}{243}\leq |T_p|\leq\frac{|G_p|}{2}$ and $|T_pv(p)\triangle T_p|\leq\frac{C}{\sqrt{p}}|G_p|$.
Since $S_p$ is $H_p$-invariant,  $T_p\rho_p(b_j)=T_p$, for all $1\leq j\leq k-1$, and $T_p\rho_p(b_k)=T_pv(p)$. Hence,
$|T_p\rho_p(b_j)\triangle T_p|\leq\frac{C}{\sqrt{p}}|G_p|$, for every $1\leq j\leq k$.

Finally, since $\ker(\rho_p)\subset\ker(\xi_p)\cap\Lambda$, the moreover assertion follows.
 \hfill$\blacksquare$


\subsection{Construction of asymptotic homomorphisms} Assume the notation from \ref{G_p} and \ref{setting}, and let $\varphi_p:\Gamma\rightarrow G_p$ and  $\rho_p:\Lambda\rightarrow G_p$ be the homomorphisms provided by Lemmas \ref{ALW01} and \ref{ALW01II}.

\begin{lemma}\label{asymptotic}  Let $t=\pm 1$ be a generator of $\mathbb Z$. Then there exists an asymptotic homomorphism $\sigma_p:(\Gamma*\mathbb Z)\times\Lambda\rightarrow\emph{Sym}(G_p)$, where $p$ is a large enough prime with $p\equiv 1\pmod{3}$, so that 
\begin{enumerate}
\item $\sigma_p(g,h)x=\varphi_p(g)x\rho_p(h)^{-1}$, for every $g\in\Gamma, h\in\Lambda, x\in G_p$,
\item $|\{x\in G_p\mid\sigma_p(t,e)x=x\}|\geq\frac{|G_p|}{3}$,
\item $\max\{\emph{d}_{\emph{H}}(\sigma_p(t,e)\circ\sigma_p(e,h),\sigma_p(e,h)\circ\sigma_p(t,e))\mid h\in\Lambda\}\geq \frac{1}{243}$, and
\item $|\sigma_p(t,e)(A_ph)\triangle A_ph'|\geq\frac{1}{243}|A_p|$, for every $h,h'\in H_p$.
\end{enumerate}
\end{lemma}

{\it Proof.} Define a homomorphism $\sigma_p:\Gamma\times\Lambda\rightarrow\text{Sym}(G_p)$ using the formula from (1). 
In order to extend $\sigma_p$ to an asymptotic homomorphism $\sigma_p:(\Gamma*\mathbb Z)\times\Lambda\rightarrow\text{Sym}(G_p)$ we will define $\sigma_p(t,e)\in\text{Sym}(G_p)$ such that $\lim\limits_{p\rightarrow\infty}\text{d}_{\text{H}}(\sigma_p(t,e)\sigma_p(e,h),\sigma_p(e,h)\sigma_p(t,e))=0$, for any $h\in\Lambda$.

To this end, let $a_p=(0,0,1...,1)\in F_3^{p+1}$. Since $3\mid p-1$, we get that $a_p\in A_p$.
If $x\in F_3^{p+1}$ and $\tilde x=x-a(p)$, then $n_1(\tilde x)-n_2(\tilde x)\leq n_0(x)-n_1(x)+4$.
This observation together with the definition \eqref{S_p} of $S_p$ implies that $(a_p+S_p)\cap S_p=\emptyset$. Since $T_p=S_p\cdot H_p$, we get that $a_pT_p\cap T_p=\emptyset$. 

Let $h_p\in H_p$ such that $h_p^2\not=e$. Since $h_pT_p=T_p$, we get that $a_ph_pT_p=a_pT_p$.
We can therefore define $\sigma_p(t,e)\in\text{Sym}(G_p)$ by letting $$\sigma_p(t,e)x=\begin{cases}\text{$a_ph_px$, if $x\in T_p$}, \\ \text{$(a_ph_p)^{-1}x$, if $x\in a_pT_p$,}\\ \text{$x$, otherwise.}  \end{cases}$$
Equivalently, for every $a\in A_p$ and $h\in H_p$ we have that \begin{equation}\label{sigma_p}\sigma_p(t,e)(ah)=\begin{cases}\text{$a_ph_pah$, if $a\in S_p$}, \\ \text{$(a_ph_p)^{-1}(ah)$, if $a\in A_pS_p$}, \\ \text{$ah$, otherwise} \end{cases}\end{equation}

By Lemma \ref{ALW01II} we have that $\lim\limits_{p\rightarrow\infty}\frac{|T_p\rho_p(h)\triangle T_p|}{|G_p|}=0$, for every $h\in\Lambda$. By arguing as in the proof of \cite[Lemma 6.1]{Io19b} we get that $\lim\limits_{p\rightarrow\infty}\text{d}_{\text{H}}(\sigma_p(t,e)\sigma_p(e,h),\sigma_p(e,h)\sigma_p(t,e))=0$, for any $h\in\Lambda$. Hence, $(\sigma_p)$ defines an asymptotic homomorphism which satisfies condition (1) by construction. 

Since $|T_p|\leq\frac{|G_p|}{3}$ by Lemma \eqref{ALW01II}, condition  (2) is satisfied.
 Moreover, since $(a_ph_p)^2\not=e$, formula (6.3) from the proof of \cite[Lemma 6.1]{Io19b} implies that \begin{equation}\label{lower}\text{$\text{d}_{\text{H}}(\sigma_p(t,e)\sigma_p(e,h),\sigma_p(e,h)\sigma_p(t,e))\geq \frac{2|T_p\setminus T_p\rho_p(h)|}{|G_p|}$, for every $h\in\Lambda$.}\end{equation}
Since $\sum_{k\in G_p}|T_p\setminus T_pk|=\sum_{k\in G_p}(|T_p|-|T_p\cap T_pk|)=(|G_p|-|T_p|)\cdot |T_p|$, there is $k\in G_p$ such that $|T_p\setminus T_pk|\geq (1-\frac{|T_p|}{|G_p|})\cdot\frac{|T_p|}{|G_p|}\cdot |G_p|$. By Lemma \ref{ALW01II} we get that $|T_p\setminus T_pk|\geq \frac{1}{2}\cdot\frac{1}{243}\cdot |G_p|.$ Since $\rho_p$ is onto,  there is $h\in \Lambda$ such that $k=\rho_p(h)$ and hence $|T_p\setminus T_p\rho_p(h)|\geq\frac{1}{2}\cdot\frac{1}{243}\cdot |G_p|$. In combination with \eqref{lower}, we deduce condition (3).

To prove condition (4), fix $h,h'\in H$. Since $h_p^2\not=e$, we have $h_ph\not=h'$ or $h_p^{-1}h\not=h'$. If $h_ph\not=h'$, then \eqref{sigma_p} implies that $S_p\subset \{a\in A_p\mid \sigma_p(t,e)(ah)\notin A_ph'\}$. If $h_p^{-1}h\not=h'$, then \eqref{sigma_p} implies that $a_pS_p\subset\{a\in A_p\mid\sigma_p(t,e)(ah)\notin A_ph'\}$. In either case, $|\{a\in A_p\mid \sigma_p(t,e)(ah)\notin A_ph'\}|\geq |S_p|$ and thus $|\sigma_p(t,e)(A_ph)\setminus A_ph'|\geq |S_p|$. By equation \eqref{bounds}, this implies condition (4).
 \hfill$\blacksquare$
 
Note that condition (2) of Lemma \ref{asymptotic} implies that $(\sigma_p)$ is not a sofic approximation of $(\Gamma*\mathbb Z)\times\Lambda$. In the next section, we will first build a sofic approximation $(\widetilde{\sigma}_p)$ of $(\Gamma*\mathbb Z)\times\Lambda$ out of $(\sigma_p)$,  which we will use to show that the conclusion of our main theorem holds for $(\Gamma*\mathbb Z)\times\Lambda=\mathbb F_m\times\mathbb F_k$.

\section{Proof of the main theorem} \label{4}
This section is devoted to the proof of the main theorem.
 We will first prove the conclusion when $m\geq 5$ and $k\geq 3$. To this end, put $\Gamma=\mathbb F_{m-1}, \Sigma=\mathbb F_m$ and $\Lambda=\mathbb F_k$. We assume the notation from \eqref{G_p} and \eqref{setting}:  $H_p=\text{PSL}_2(F_p)$, $G_p=A_p\rtimes H_p$, $K_p=\text{PSL}_2(F_{r_p})$, $\widetilde G_p=G_p\times K_p$, $\psi_p:\text{PSL}_2(\mathbb Z)\rightarrow K_p$ is the quotient homomorphism, where $p<r_p$ are primes and $p\equiv 1\pmod{3}$.
 
 {\it In the first part of the proof}, we construct a sofic approximation  $\widetilde\sigma_p:\Sigma\times\Lambda\rightarrow\text{Sym}(\widetilde G_p)$ of $\Sigma\times\Lambda$.
 View $\Sigma$, and thus $\Gamma$, as a subgroup of $\text{PSL}_2(\mathbb Z)$.
 Let $\varphi_p:\Gamma\rightarrow G_p$, $\widetilde\varphi_p:\Gamma\rightarrow\widetilde G_p$ and  $\rho_p:\Lambda\rightarrow G_p$ be the onto homomorphisms given by Lemmas \ref{ALW01} and \ref{ALW01II}. 
Recall that $\widetilde\varphi_p(g)=(\varphi_p(g),\psi_p(g))$, for every $g\in\Gamma$. 
Let $\sigma_p:\Sigma\times\Lambda=(\Gamma*\mathbb Z)\times\Lambda\rightarrow\text{Sym}(G_p)$ be the asymptotic homomorphism provided by Lemma \ref{asymptotic}.  As  therein, we denote by $t=\pm 1$ a generator of $\mathbb Z$.

 Next, for any large enough prime $p$ with $p \equiv 1\pmod{3}$, we  define onto homomorphisms $\zeta_p:\Lambda\rightarrow K_p$ and $\widetilde\rho_p:\Lambda\rightarrow \widetilde G_p$. Fix a decomposition $\Lambda=\Delta*\mathbb Z$, where $\Delta=\mathbb F_{k-1}$ and view $\Delta$ as a subgroup of $\text{PSL}_2(\mathbb Z)$. 
  Define $\zeta_p:\Lambda\rightarrow K_p$ by $\zeta_p(h)=\psi_p(h)$, if $h\in\Delta$, and $\zeta_p(h)=e$, if $h\in\mathbb Z$.
Define $\widetilde\rho_p:\Lambda\rightarrow\widetilde G_p$ by  $\widetilde\rho_p(h)=(\rho_p(h),\zeta_p(h))$, for all  $h\in\Lambda$. 
Since $\Delta$ is non-amenable, $\zeta_p$ is onto, for large enough $p$. Since $G_p$ and $K_p$ are simple non-isomorphic groups and $\rho_p$ is onto, Goursat's lemma implies that $\widetilde\rho_p$ is onto, for large enough $p$.

We are now ready to define
 $\widetilde\sigma_p:\Sigma\times\Lambda\rightarrow\text{Sym}(\widetilde G_p)$ by letting for $g\in\Sigma,h\in\Lambda, x\in G_p$ and $y\in K_p$ 
\begin{equation}\label{defp}
\widetilde{\sigma}_p(g,h)(x,y)=(\sigma_p(g,h)x,\psi_p(g)y\zeta_p(h)^{-1})
\end{equation} 
Since  Lemma \ref{asymptotic} gives that  $\sigma_p(g,h)x=\varphi_p(g)x\rho_p(h)^{-1}$, for all $g\in\Gamma, h\in\Lambda,x\in G_p$, we derive that \begin{equation}\label{tres} \text{$\widetilde\sigma_p(g,h)x=\widetilde\varphi_p(g)x\widetilde\rho_p(h)^{-1}$, for all $g\in\Gamma,h\in\Lambda,x\in\widetilde G_p$}. \end{equation}

\begin{claim}\label{soficapp} $(\widetilde\sigma_p)$ is a sofic approximation of $\Sigma\times\Lambda$.
\end{claim}

{\it Proof of Claim \ref{soficapp}.} It is clear that $(\widetilde\sigma_p)$ is an asymptotic homomorphism of $\Sigma\times\Lambda$.
To see that it is a sofic approximation, let $C_L(g)=\{y\in L\mid gy=yg\}$ be the centralizer of an element $y$ of a group $L$. For every prime $q$ and $g\in \text{SL}_2(F_q)\setminus\{\pm I\}$, we have $|C_{\text{SL}_2(F_q)}(g)|\leq \frac{|\text{SL}_2(F_q)|}{q-1}$. This implies that \begin{equation}\label{sofique}\text{$|C_{\text{PSL}_2(F_q)}(g)|\leq \frac{|\text{PSL}_2(F_q)|}{2(q-1)}$, for every $g\in\text{PSL}_2(F_q)\setminus\{e\}$.}\end{equation}
On the other hand,  for all $(g,h)\in\Sigma\times\Lambda$ we have\begin{equation}\label{sofiques}
\frac{|\{x\in\widetilde G_p\mid \widetilde\sigma_p(g,h)x=x\}|}{|\widetilde G_p|}\leq\frac{|\{y\in K_p\mid\psi_p(g)y\zeta_p(h)^{-1}=y\}|}{|K_p|}\leq\frac{|C_{K_p}(\psi_p(g))|}{|K_p|}.
\end{equation}
If $g\not=e$, then $\psi_p(g)\not=e$, for large enough  $p$. Therefore, by combining \eqref{sofique} and \eqref{sofiques} we get that $\frac{|\{x\in\widetilde G_p\mid \widetilde\sigma_p(g,h)x=x\}|}{|\widetilde G_p|}\leq\frac{1}{2(r_p-1)}$, for large enough $p$. Thus, $\lim\limits_{p\rightarrow\infty}\text{d}_{\text{H}}(\widetilde\sigma_p(g,h),\text{Id}_{\widetilde G_p})=1$, for all $(g,h)\in\Sigma\times\Lambda$ with $g\not=e$. If $h\in\Lambda\setminus\{e\}$, then $\widetilde\sigma_p(e,h)(x,y)=(x\rho_p(h)^{-1},y\zeta_p(h)^{-1})$. 
Lemma \ref{ALW01II} gives that $\rho_p(h)\not=e$ and thus $\text{d}_{\text{H}}(\widetilde\sigma_p(e,h),\text{Id}_{\widetilde G_p})=1$, for  large enough $p$. This altogether proves that  $(\widetilde\sigma_p)$ is a sofic approximation of $\Sigma\times\Lambda$. \hfill$\square$

{\it In the rest of the proof}, we will show that the sofic approximation $(\widetilde\sigma_p)$ of $\Sigma\times\Lambda$ satisfies the conclusion of the main theorem. Towards this goal, let $\tau_p:\Sigma\times\Lambda\rightarrow\text{Sym}(Y_p)$ be a sequence homomorphisms, for some finite sets $Y_p$, for which there exist maps $\theta_p:\widetilde G_p\rightarrow Y_p$ such that 
\begin{enumerate}[label=(\roman*)]
\item\label{uno} $\text{d}_{\text H}(\theta_p\circ\widetilde\sigma_p(g),\tau_p(g)\circ\theta_p)\rightarrow 0$, for every $g\in\Sigma\times\Lambda$, and 
\item\label{dos} $\text{d}_{\text{H}}(\theta_p\circ\widetilde\sigma_p(g),\theta_p)\rightarrow 1$, for every $g\in(\Sigma\times\Lambda)\setminus\{e\}$.
\end{enumerate}

 For every $p$ and $y\in Y_p$, we denote $X_p^y=\theta_p^{-1}(\{y\})$. We continue with the following:
 
 \begin{claim}\label{normal}
 For every $p$  there exist  a normal subgroup $N_p<\widetilde G_p$, a subset $S_p\subset Y_p$ and a map $\omega_p:S_p\rightarrow\widetilde G_p$ such that \begin{equation}\label{N_p}\text{$\lim\limits_{p\rightarrow\infty}\frac{1}{|\widetilde G_p|}\cdot\sum_{y\in S_p}|X_p^y\triangle \omega_p(y)N_p|=0$ \;\;\; and}\end{equation} 
\begin{equation}\label{N_pp}\text{$\lim\limits_{p\rightarrow\infty}\frac{1}{|\widetilde G_p|}\cdot\sum_{y\notin S_p}|X_p^y|=0$.}\end{equation}
 Moroever, for every $g\in\Sigma\times\Lambda$, we have that $\displaystyle{\lim\limits_{p\rightarrow\infty}\frac{|N_p|\cdot |S_p\cap \tau_p(g)^{-1}S_p|}{|\widetilde G_p|}=1}$.
 \end{claim}

{\it Proof of Claim  \ref{normal}}. If $g\in\Sigma\times\Lambda$, then $$\bigcup_{y\in Y_p}\big(\widetilde\sigma_p(g)X_p^y\triangle X_p^{\tau_p(g)(y)}\big)\subset\{x\in \widetilde G_p\mid\theta_p(\widetilde\sigma_p(g)^{-1}x)\not=\tau_p(g)^{-1}(\theta_p(x))\}.$$

Since $\text{d}_{\text H}(\theta_p\circ\widetilde\sigma_p(g)^{-1},\tau_p(g)^{-1}\circ\theta_p)\rightarrow 0$ by \ref{uno}, we deduce that
the partition $\{X_p^y\}_{y\in Y_p}$ of $\widetilde G_p$ is almost invariant under $\widetilde\sigma_p$, in the following sense:  \begin{equation}\label{ainv}\text{$\lim\limits_{p\rightarrow\infty}\frac{1}{|\widetilde G_p|}\cdot \sum_{y\in Y_p}|\widetilde\sigma_p(g)X_p^y\triangle X_p^{\tau_p(g)(y)}|=0$, for every $g\in\Sigma\times\Lambda$.}
\end{equation}

By Lemma \ref{ALW01}, $\widetilde\varphi_p:\Gamma\rightarrow\widetilde G_p$ is an onto homomorphism such that $\Gamma$ has property $(\tau)$ with respect to $\{\ker(\widetilde\varphi_p)\}$. Moreover, \eqref{tres} gives that
  $\widetilde\sigma_p(g,e)(x)=\widetilde\varphi_p(g)x$, for all $g\in\Gamma$ and $x\in \widetilde G_p$. 
Since equation \eqref{ainv} holds for all $g\in\Gamma$, we can apply Lemma \ref{Io18} to deduce the existence of a subgroup $N_p<\widetilde G_p$, a subset $S_p\subset Y_p$ and a map $\omega_p:S_p\rightarrow\widetilde G_p$, for every $p$, such that \eqref{N_p} and \eqref{N_pp} hold.

To finish the proof of the claim, it remains to show that $N_p<\widetilde G_p$ is a normal subgroup and the moreover assertion holds. To this end,  combining \eqref{N_pp} and \eqref{ainv} gives that$$\text{$\lim\limits_{p\rightarrow\infty}\frac{1}{|\widetilde G_p|}\cdot\sum_{y\notin\tau_p(g)^{-1}S_p}|X_p^y|=0$, for every $g\in\Sigma\times\Lambda$.}$$ Since $\sum_{y\in Y_p}|X_p^y|=|\widetilde G_p|$, this together with \eqref{N_pp} implies that
\begin{equation}\label{N_ppp} 
\text{$\lim\limits_{p\rightarrow\infty}\frac{1}{|\widetilde G_p|}\cdot\sum_{y\in S_p\cap\tau_p(g)^{-1}S_p}|X_p^y|=1$, for every $g\in\Sigma\times\Lambda$.}
\end{equation}
By combining \eqref{N_p} and \eqref{N_ppp}, the moreover assertion follows.

On the other hand, combining \eqref{N_p} and \eqref{ainv} gives that \begin{equation}\label{N_pppp}\text{$\lim\limits_{p\rightarrow\infty}\frac{1}{|\widetilde G_p|}\cdot\sum_{y\in S_p\cap\tau_p(g)^{-1}S_p}|\widetilde\sigma_p(g)(\omega_p(y)N_p)\triangle\omega_p(\tau_p(g)(y))N_p|=0$, for every $g\in\Sigma\times\Lambda$.} \end{equation}

Let $h\in\Lambda$. Then  \eqref{N_pppp} and the moreover assertion imply that for every large enough $p$, there is $y\in S_p\cap\tau_p(e,h)^{-1}S_p$ such that $ |\widetilde\sigma_p(e,h)(\omega_p(y)N_p)\triangle\omega_p(\tau_p(e,h)(y))N_p|<|N_p|$.
 Since by \eqref{tres} we have $\widetilde\sigma_p(e,h)(x)=x\widetilde\rho_p(h)^{-1}$, for all $x\in\widetilde G_p$, we get
 $|\omega_p(y)N_p\widetilde\rho_p(h)^{-1}\triangle \omega_p(\tau_p(g)(y))N_p|<|N_p|$. Thus, if we put $a=\widetilde\rho_p(h)\omega_p(y)^{-1}\omega_p(\tau_p(g)(y))\in\widetilde G_p$, then
$|\widetilde\rho_p(h)N_p\widetilde\rho_p(h)^{-1}\triangle aN_p|<|N_p|$. This further implies that $|\widetilde\rho_p(h)N_p\widetilde\rho_p(h)^{-1}\cap N_p|>\frac{|N_p|}{2}$. Since $N_p$ is a finite group, it follows that $\widetilde\rho_p(h)N_p\widetilde\rho_p(h)^{-1}=N_p$, for every large enough prime $p$. Since this holds for every $h\in\Lambda$, $\Lambda$ is finitely generated and $\widetilde\rho_p:\Lambda\rightarrow\widetilde G_p$ is onto, we derive that $N_p<\widetilde G_p$ is a normal subgroup. \hfill$\square$

Since $\widetilde G_p=G_p\times K_p$ and $K_p$ is a simple group, Lemma \ref{irred}  implies that $N_p$ is  one of the following six groups: $N_p^1=\{e\}, N_p^2=\{e\}\times K_p, N_p^3=A_p\times K_p, N_p^4=\widetilde G_p, N_p^5=G_p\times\{e\}$ or $N_p^6=A_p\times\{e\}$, for every large enough prime $p$. 
We continue with the following:

\begin{claim}\label{trivial}
$N_p=\{e\}$, for every large enough prime $p$.
\end{claim} 

{\it Proof of Claim \ref{trivial}.} Assume that the claim is false. Then, after replacing $(\widetilde\sigma_p)$ with a subsequence, we may assume that there is $2\leq i\leq 6$ such that $N_p=N_p^i$, for every $p$. We will prove that each of these five possibilities leads to a contradiction. To this end, note that if $g\in \Sigma\times\Lambda$, then $\cup_{y\in Y_p}(X_p^y\cap\widetilde\sigma_p(g)X_p^y)\subset\{x\in\widetilde G_p\mid(\theta_p\circ\widetilde\sigma_p(g)^{-1})(x)=\theta_p(x)\}$. Thus, by using \ref{dos}, we get that 
\begin{equation}
\label{sofic}\text{$\lim\limits_{p\rightarrow\infty}\frac{1}{|\widetilde G_p|}\cdot\sum_{y\in Y_p}|X_p^y\cap\widetilde\sigma_p(g)X_p^y|=0$, for every $g\in (\Sigma\times\Lambda)\setminus\{e\}$.}
\end{equation}

\begin{case} $N_p=L_p\times K_p$, where $L_p$ is one of the groups $\{e\}, A_p$ or $G_p$.
\end{case}

In this case, we will derive a contradiction by using condition (2) of Lemma \ref{asymptotic}.
Note first that for every $y\in Y_p$, we have that $\omega_p(y)N_p=V_p^y\times K_p$, for some set $V_p^y\subset G_p$. 
By \eqref{defp} we have that $\widetilde\sigma_p(t,e)(x,y)=(\sigma_p(t,e)x,\psi_p(t)y)$, for all $x\in G_p$ and $y\in K_p$. Thus, we get that $$\widetilde\sigma_p(t,e)(\omega_p(y)N_p)=\sigma_p(t,e)V_p^y\times K_p=(\sigma_p(t,e)\times\text{Id}_{K_p})(\omega_p(y)N_p).$$ This implies that $\omega_p(y)N_p\cap\widetilde\sigma_p(t,e)(\omega_p(y)N_p)\supset \{x\in\omega_p(y)N_p\mid  (\sigma_p(t,e)\times\text{Id}_{K_p})x=x\}$ and hence \begin{align*}\sum_{y\in S_p}|\omega_p(y)N_p\cap\widetilde\sigma_p(t,e)(\omega_p(y)N_p)|&\geq|\{x\in\bigcup_{y\in S_p}\omega_p(y)N_p\mid (\sigma_p(t,e)\times\text{Id}_{K_p})(x)=x\}|\\ &\geq \{x\in\widetilde G_p\mid(\sigma_p(t,e)\times\text{Id}_{K_p})(x)=x\}|-{|\widetilde G_p\setminus(\bigcup_{y\in S_p}\omega_p(y)N_p))|}.
\end{align*}
Since $|\{x\in G_p\mid\sigma_p(t,e)x=x\}|\geq \frac{|G_p|}{3}$ by condition (2) of Lemma \ref{asymptotic},  while \eqref{N_p} and \eqref{N_pp} imply that $\lim\limits_{p\rightarrow\infty}\frac{|\widetilde G_p\setminus(\cup_{y\in S_p}\omega_p(y)N_p)|}{|\widetilde G_p|}=0$, we conclude that $$\limsup_{p\rightarrow\infty}\frac{1}{|\widetilde G_p|}\cdot \sum_{y\in S_p}|\omega_p(y)N_p\cap\widetilde\sigma_p(t,e)(\omega_p(y)N_p)|\geq \frac{1}{3}.$$ In combination with \eqref{N_p} this gives that $$\limsup_{p\rightarrow\infty}\frac{1}{|\widetilde G_p|}\cdot\sum_{y\in S_p}|X_p^y\cap\widetilde\sigma_p(t,e)X_p^y|\geq\frac{1}{3},$$ which contradicts \eqref{sofic}.

\begin{case} $N_p=G_p\times\{e\}$.
\end{case}

In this case, for every $y\in Y_p$ we have that $\omega_p(y)N_p=G_p\times W_p^y$, for some set  $W_p^y\subset K_p$. Recall that by the construction of the homomorphism $\zeta_p:\Lambda\rightarrow K_p$ we can find $h\in\Lambda\setminus\{e\}$ such that $\zeta_p(h)=e$. Then $\widetilde\rho_p(h)=(\rho_p(h),e)$ and \eqref{tres} gives that $\widetilde\sigma_g(e,h)(x,y)=(x\rho_p(h)^{-1},y)$ for all $x\in G_p$ and $y\in K_p$. Hence $\widetilde\sigma_g(e,h)(\omega_p(y)N_p)=\omega_p(y)N_p$, for every $y\in Y_p$. By reasoning as in Case 1, this implies that $\limsup_{p\rightarrow\infty}\frac{1}{|\widetilde G_p|}\cdot\sum_{y\in S_p}|X_p^y\cap\widetilde\sigma_p(e,h)X_p^y|\geq 1,$ which contradicts \eqref{sofic}.

\begin{case} $N_p=A_p\times\{e\}$.
\end{case}

In this case, 
let $x,x'\in G_p$ and $y,y'\in K_p$. Then $(x,y)N_p=A_px\times\{y\}$ and $(x',y')N_p=A_px'\times\{y'\}$. By using the definition \eqref{defp} of $\widetilde\sigma_p$ we derive that $\widetilde\sigma_p(t,e)((x,y)N_p)=\sigma_p(t,e)(A_px)\times\{\psi_p(t)y\}$. In combination with condition (4)  from Lemma \ref{asymptotic}, we conclude that $$|\widetilde\sigma_p(t,e)((x,y)N_p)\triangle (x',y')N_p|\geq  |\sigma_p(t,e)(A_px)\triangle (A_px')|\geq\frac{|A_p|}{243}=\frac{|N_p|}{243}. $$
This inequality implies that \begin{equation}\label{contra}\frac{1}{|\widetilde G_p|}\cdot\sum_{y\in S_p\cap\tau_p(g)^{-1}S_p}|\widetilde\sigma_p(t,e)(\omega_p(y)N_p)\triangle\omega_p(\tau_p(e,h)(y))N_p|\geq \frac{|N_p|\cdot |S_p\cap\tau_p(t,e)^{-1}S_p|}{243\cdot|\widetilde G_p|}.\end{equation} On the other hand, \eqref{N_pppp} and the moreover assertion of Lemma \ref{normal} imply that the left and right sides of \eqref{contra} converge to $0$ and $\frac{1}{243}$,  as $p\rightarrow\infty$, respectively. This gives a contradiction.

This altogether finishes the proof of Claim \ref{trivial}. \hfill$\square$

We next claim that $\theta_p$ is ``asymptotically one-to-one": there exists a set  $X_p\subset\widetilde G_p$ such that ${\theta_p}_{|X_p}$ is one-to-one and $\lim\limits_{p\rightarrow\infty}\frac{|X_p|}{|\widetilde G_p|}=1$.  To see this, let $S_p'=\{y\in S_p\mid |X_p^y|\not=1\}$. If $y\in S_p\setminus S_p'$, then $|X_p^y|\not=1$. Since $N_p=\{e\}$ we have that $|\omega_p(y)N_p|=1$ and thus $|X_p^y\triangle\omega_p(y)N_p|\geq\frac{|X_p^y|}{2}$. In combination with \eqref{N_p} we derive that $\lim\limits_{p\rightarrow\infty}\frac{1}{|\widetilde G_p|}\cdot\sum_{y\in S_p\setminus S_p'}|X_p^y|=0$. Together with \eqref{N_pp} this further implies that $\lim\limits_{p\rightarrow\infty}\frac{1}{|\widetilde G_p|}\cdot\sum_{y\notin S_p'}|X_p^y|=0$. Thus, if we let $X_p=\cup_{y\in S_p'}X_p^y$, then $\lim\limits_{p\rightarrow\infty}\frac{|X_p|}{|\widetilde G_p|}=1$. 
Since $|X_p^y|=1$, for every $y\in S_p'$, we also have that ${\theta_p}_{|X_p}$ is one-to-one, which proves our claim.

Finally, we will use  that $\tau_p:\Sigma\times\Lambda\rightarrow\text{Sym}(Y_p)$ is a homomorphism, for every $p$, in combination with  \cite[Theorem 5.1]{Io19b} to derive a contradiction. Consider the disjoint union $Z_p=Y_p\sqcup (\widetilde G_p\setminus X_p)$ and extend $\tau_p$ to a homomorphism $\tau_p:\Sigma\times\Lambda\rightarrow\text{Sym}(Z_p)$ by letting $\tau_p(g)_{|\widetilde G_p\setminus X_p}=\text{Id}_{\widetilde G_p\setminus X_p}$, for every $g\in\Sigma\times\Lambda$.
We also define a  one-to-one map $\Theta_p:\widetilde G_p\rightarrow Z_p$ by letting $$\Theta_p(x)=\begin{cases}\text{$\theta_p(x)$, if $x\in X_p$, and} \\ \text{$x$, if $x\in\widetilde G_p\setminus X_p$.}  \end{cases}$$
If $g\in\Sigma\times\Lambda$, then $\Theta_p(x)=\theta_p(x)$ and $\Theta_p(\widetilde\sigma_p(g)x)=\theta_p(\widetilde\sigma_p(g)x)$, for every $x\in X_p\cap\widetilde\sigma_p(g)^{-1}X_p$. From this it follows that $$\text{d}_{\text{H}}(\Theta_p\circ\widetilde\sigma_p(g),\tau_p(g)\circ\Theta_p)\leq \text{d}_{\text{H}}(\theta_p\circ\widetilde\sigma_p(g),\tau_p(g)\circ\theta_p)+\frac{|\widetilde G_p\setminus(X_p\cap\widetilde\sigma_p(g)^{-1}X_p)|}{|\widetilde G_p|}.$$
By using \ref{uno}, that $|\widetilde G_p\setminus(X_p\cap\widetilde\sigma_p(g)^{-1}X_p)|\leq 2\cdot |\widetilde G_p\setminus X_p|$ and that $\lim\limits_{p\rightarrow\infty}\frac{|X_p|}{|\widetilde G_p|}=1$, we deduce that \begin{equation}\label{embed}\text{$\lim\limits_{p\rightarrow\infty}\text{d}_{\text{H}}(\Theta_p\circ\widetilde\sigma_p(g),\tau_p(g)\circ\Theta_p)=0$, for every $g\in\Sigma\times\Lambda$}.\end{equation}

Note that $\widetilde\sigma_p(g,h)x=\widetilde\varphi_p(g)x\widetilde\rho_p(h)^{-1}$, for all $g\in\Gamma,h\in\Lambda, x\in\widetilde G_p$ by \eqref{tres}, $\widetilde\varphi_p$ is onto and $\Gamma$ has property $(\tau)$ with respect to $\{\ker(\widetilde\varphi_p)\}_p$ by Lemma \ref{ALW01}. Since $\tau_p$ is a homomorphism and $\Theta_p$ is one-to-one, for every $p$, we can apply \cite[Theorem 5.1]{Io19b} to conclude that $$\lim\limits_{p\rightarrow\infty}\Big(\max\{\text{d}_{\text{H}}(\widetilde\sigma_p(t,e)\circ\widetilde\sigma_p(e,h),\widetilde\sigma_p(e,h)\circ\widetilde\sigma_p(t,e))\mid h\in\Lambda\}\Big)=0.$$ Using the definition \ref{defp} of $\widetilde\sigma_p$, this can be equivalently written as  $$\lim\limits_{p\rightarrow\infty}\Big(\max\{\text{d}_{\text{H}}(\sigma_p(t,e)\circ\sigma_p(e,h),\sigma_p(e,h)\circ\sigma_p(t,e))\mid h\in\Lambda\}\Big)=0,$$ which contradicts condition (3) from Lemma \ref{asymptotic}. 

This finishes the proof of  our main theorem in the case $m\geq 5$ and $k\geq 3$. 
In the general case, when $m,k\geq 2$, note first that $\mathbb F_{m+3}\times\mathbb F_{k+1}$ satisfies the conclusion by the above. Since $\mathbb F_{m+3}\times\mathbb F_{k+1}$ can be realized as a finite index subgroup of $\mathbb F_m\times\mathbb F_k$,  the following lemma implies that the conclusion also holds for $\mathbb F_m\times\mathbb F_k$. \hfill$\blacksquare$

Consider the following property for a sofic approximation $\sigma_n:\Gamma_0\rightarrow \text{Sym}(X_n)$ of a countable group $\Gamma_0$: ($\diamond)$  there are a sofic approximation $\tau_n:\Gamma_0\rightarrow\text{Sym}(Y_n)$ and maps $\theta_n:X_n\rightarrow Y_n$ such that
\begin{enumerate}[label={(\alph*)}]
\item $\tau_n$ is a homomorphism, for every $n\in\mathbb N$,
\item $\text{d}_{\text{H}}(\theta_n\circ\sigma_n(g),\tau_n(g)\circ\theta_n)\rightarrow 0$, for every $g\in\Gamma_0$, and
\item $\text{d}_{\text{H}}(\theta_n\circ\sigma_n(g),\theta_n)\rightarrow 1$, for every $g\in\Gamma_0\setminus\{e\}$.
\end{enumerate}

Let $\Gamma_0<\Gamma$ be a finite index inclusion of countable groups. 
The next lemma shows that if $\Gamma_0$ has a sofic approximation which fails $(\diamond)$, then so does  $\Gamma$. The proof of this fact relies on an induction argument, following closely \cite[Section 3.3]{Io19b}.
Let $s:\Gamma/\Gamma_0\rightarrow\Gamma$ be a map such that $s(e\Gamma_0)=e$ and $s(g\Gamma_0)\in g\Gamma_0$,  for every $g\in \Gamma$. Then $c:\Gamma\times\Gamma/\Gamma_0\rightarrow\Gamma_0$ given by $c(g,h\Gamma_0)=s(gh\Gamma_0)^{-1}g s(h\Gamma_0)$ is a cocycle for the left multiplication action $\Gamma\curvearrowright \Gamma/\Gamma_0$.

Let $\sigma_n:\Gamma_0\rightarrow\text{Sym}(X_n)$ be a sofic approximation of $\Gamma_0$ and define $\text{Ind}_{\Gamma_0}^{\Gamma}(\sigma_n):\Gamma\rightarrow\text{Sym}(\Gamma/\Gamma_0\times X_n)$ by letting $\text{Ind}_{\Gamma_0}^{\Gamma}(\sigma_n)(g)(h\Gamma_0,x)=(gh\Gamma_0,\sigma_n(c(g,h\Gamma_0))x)$, for every $g\in\Gamma,h\Gamma_0\in\Gamma/\Gamma_0$ and $x\in X_n$. 
The proof of \cite[Lemma 3.3]{Io19b} shows that $\text{Ind}_{\Gamma_0}^{\Gamma}(\sigma_n)$ is a sofic approximation of $\Gamma$. 

\begin{lemma}\label{fin}
 If the induced sofic approximation $(\emph{Ind}_{\Gamma_0}^{\Gamma}(\sigma_n))$  satisfies $(\diamond$), then $(\sigma_n)$ satisfies $(\diamond)$.
\end{lemma}

{\it Proof.}  Assume that $\widetilde\sigma_n:=\text{Ind}_{\Gamma_0}^{\Gamma}(\sigma_n):\Gamma\rightarrow\text{Sym}(\widetilde X_n)$ satisfies $(\diamond)$, where  $\widetilde X_n=\Gamma/\Gamma_0\times X_n$. Let $\tau_n:\Gamma\rightarrow\text{Sym}(Y_n)$ be a sofic approximation by homomorphisms and $\widetilde\theta_n:\widetilde X_n\rightarrow Y_n$ be maps such that $\text{d}_{\text{H}}(\widetilde\theta_n\circ\widetilde\sigma_n(g),\tau_n(g)\circ\widetilde\theta_n)\rightarrow 0$, for every $g\in\Gamma$, and
 $\text{d}_{\text{H}}(\widetilde\theta_n\circ\widetilde\sigma_n(g),\widetilde\theta_n)\rightarrow 1$, for every $g\in\Gamma\setminus\{e\}$. If $g\in\Gamma_0$, then $\widetilde\sigma_n(g)$ leaves $e\Gamma_0\times X_n$ invariant and  $\widetilde\sigma_n(g)(e\Gamma_0,x)=(e\Gamma_0,\sigma_n(g)x)$, for every $x\in X_n$. Thus, the restriction of ${\widetilde{\sigma_n}}_{|\Gamma_0}$ 
 to $e\Gamma_0\times X_n$ can be identified to $\sigma_n$. Denote by $\theta_n:e\Gamma_0\times X_n\rightarrow Y_n$ the restriction of $\widetilde\theta_n$ to $e\Gamma_0\times X_n$. Then it is clear that for every $g\in\Gamma_0$ we have that $$\text{$\text{d}_{\text{H}}(\theta_n\circ\widetilde{\sigma_n}(g)_{|e\Gamma_0\times X_n},\tau_n(g)\circ\theta_n)\leq [\Gamma:\Gamma_0]\cdot \text{d}_{\text{H}}(\widetilde\theta_n\circ\widetilde\sigma_n(g),\tau_n(g)\circ\widetilde\theta_n)$ and}$$ $$1-\text{$\text{d}_{\text{H}}(\theta_n\circ\widetilde{\sigma_n}(g)_{|e\Gamma_0\times X_n},\theta_n)\leq [\Gamma:\Gamma_0]\cdot \big(1-\text{d}_{\text{H}}(\widetilde\theta_n\circ\widetilde\sigma_n(g),\widetilde\theta_n)\big)$.}$$
These inequalities imply that the maps $\theta_n$ witness that $(\sigma_n)$ satisfies $(\diamond$).
\hfill$\blacksquare$

\end{document}